\documentclass[reqno]{amsart}
\usepackage{amssymb}

\theoremstyle{remark}
\newtheorem{remark}{Remark}

\numberwithin{equation}{section}

\begin{document}

\title[Branching Rules]
{Branching Rules for  \\ Symmetric Hypergeometric Polynomials}

\author{J.F.  van Diejen}

\address{
Instituto de Matem\'atica y F\'{\i}sica, Universidad de Talca,
Casilla 747, Talca, Chile}

\email{diejen@inst-mat.utalca.cl}

\author{E. Emsiz}

\address{
Facultad de Matem\'aticas, Pontificia Universidad Cat\'olica de Chile,
Casilla 306, Correo 22, Santiago, Chile}
\email{eemsiz@mat.puc.cl}


\thanks{This work was supported in part by the {\em Fondo Nacional de Desarrollo
Cient\'{\i}fico y Tecnol\'ogico (FONDECYT)} Grants \# 1130226 and  \# 1141114.}

\date{October 2015}

\begin{abstract}
Starting from a recently found branching formula for the six-parameter family of symmetric Macdonald-Koornwinder polynomials, we arrive by degeneration at corresponding branching rules for symmetric hypergeometric orthogonal polynomials of Wilson, continuous Hahn, Jacobi, Laguerre, and Hermite type. 
\end{abstract}

\maketitle

\section{Introduction}\label{sec1}
Let $M_\lambda(x_1,\ldots,x_n)$, $
\lambda\in\Lambda_n:=\{(\lambda_1,\ldots, \lambda_n)\in\mathbb{Z}^n\mid \lambda_1\geq\lambda_2\geq \cdots \geq\lambda_n\geq 0\} $
denote the monomial basis of an algebra of symmetric (trigonometric) polynomials. 
We are concerned with branching rules for families of hypergeometric polynomials of the form
\begin{subequations}
\begin{equation}\label{pa}
P_\lambda (x_1,\ldots ,x_n)=M_\lambda (x_1,\ldots ,x_n)+\sum_{\substack{\mu\in\Lambda_n \\ \mu< \lambda}} c_{\lambda ,\mu} M_\mu(x_1,\ldots ,x_n)\qquad (\lambda\in\Lambda_n),
\end{equation}
with $c_{\lambda ,\mu} \in\mathbb{C}$ such that
\begin{equation}\label{pb}
\int_{\mathcal{D}} P_\lambda (x_1,\ldots,x_n) \overline{M_\mu (x_1,\ldots,x_n)}\Delta (x_1,\ldots ,x_n) \text{d}x_1\cdots \text{d}x_n=0 \quad\text{if}\quad \mu<\lambda .
\end{equation}
\end{subequations}
Here the partitions are partially ordered in accordance with the (nonhomogeneous) dominance order
$$
\mu\leq \lambda \quad\text{iff}\quad \mu_1+\cdots+\mu_k\leq\lambda_1+\cdots+\lambda_k \quad\text{for}\ k=1,\ldots,n,
$$
and the polynomial family is characterized by
a Selberg-type orthogonality weight function $\Delta (x_1,\ldots ,x_n)$  supported on   $\mathcal{D}= [-\pi ,\pi]^n$ (in the case of trigonometric polynomials) or $\mathcal{D}= \mathbb{R}^n$ (otherwise).

The branching rules under consideration are expansion formulas for the polynomials in $n+1$ variables in terms of the $n$-variable polynomials
of the form
\begin{equation}\label{bf}
P_\lambda (x_1,\ldots ,x_n,x)=
\sum_{\substack{\mu\in\Lambda_n\\ \mu \preceq \lambda}}  P_\mu (x_1,\ldots ,x_n )
P_{\lambda/\mu}(x)\qquad (\lambda\in\Lambda_{n+1}),
\end{equation}
where $\mu \preceq \lambda$ iff there exists a $\nu\in\Lambda_n$ with $\mu\subset\nu\subset\lambda$ such that the skew diagrams $\lambda/\nu$ and $\nu/\mu$ are horizontal strips. Here $\Lambda_n$ is thought of as being embedded in $\Lambda_{n+1}$  `by adding a part of size zero', and
we recall
that for $\lambda,\mu\in\Lambda_n$ one has that $\mu\subset \lambda$ iff $\mu_j\leq \lambda_j$ for $j=1,\ldots,n$, while the corresponding skew diagram $\lambda/\mu$ is a horizontal strip provided the parts of $\lambda$ and $\mu$ interlace as follows:
$$
\lambda_1\geq\mu_1\geq\lambda_2\geq\mu_2\geq\cdots \geq\lambda_n\geq\mu_n .
$$
By iterating the branching rule \eqref{bf}
\begin{align}
&P_\lambda (x_1,\ldots ,x_n,x_{n+1},\ldots, x_{n+l })= \\
&\sum_{\substack{\mu^{(n+i)}\in\Lambda_{n+i}, i=0,\ldots ,l\\ \mu^{(n)} \preceq \mu^{(n+1)}\preceq \cdots \preceq \mu^{(n+l)}=\lambda}}  P_{\mu^{(n)}} (x_1,\ldots ,x_n )
\prod_{1\leq i\leq l} P_{\mu^{(n+i)}/\mu^{(n+i-1)}}(x_{n+i}) ,\nonumber
\end{align}
one can build the polynomials in $n$ variables
\begin{equation}\label{tf}
P_\lambda (x_1,\ldots ,x_n)= 
\sum_{\substack{\mu^{(i)}\in\Lambda_{i}, i=1,\ldots ,n\\  \mu^{(1)}\preceq \mu^{(2)}\preceq \cdots \preceq \mu^{(n)}=\lambda}} P_{\mu^{(1)}}(x_1) \prod_{1< i\leq n} P_{\mu^{(i)}/\mu^{(i-1)}}(x_i)
\end{equation}
starting from the corresponding one-variable polynomials
$P_{m}(x) $, $m=0,1,2,\ldots$.

To make the above construction effective, it suffices to determine the branching polynomials
$P_{\lambda/\mu}(x)$ in Eq. \eqref{bf} explicitly. In \cite{die-ems:branching} this was realized for the symmetric
Macdonald-Koornwinder polynomials
\cite{koo:askey-wilson,mac:orthogonal}, which for $n=1$ amount to the Askey-Wilson polynomials \cite{koe-les-swa:hypergeometric}. 
Our present aim is to degenerate from the Askey-Wilson level and identify the pertinent branching polynomials needed for the explicit construction of
the symmetric hypergeometric polynomials of 
Wilson type \cite{die:multivariable,die:properties,zha:spherical,gro:multivariable}, continuous Hahn type \cite{die:multivariable,die:properties},
Jacobi type \cite{vre:formulas,hec-opd:hypergeometric,hec:hypergeometric,mac:hypergeometric,deb:systeme,las:jacobi,bee-opd:certain,bak-for:calogero-sutherland,die-lap-mor:determinantal,zha:spherical,oko-ols:limits,dum-ede-shu:MOPS,ser-ves:bc,hal-lan:unified,koo:okounkov}, Laguerre type \cite{mac:hypergeometric,las:laguerre,bak-for:calogero-sutherland,die:confluent,xu:orthogonal,ari-dav-ola:differential,dum-ede-shu:MOPS,hal-lan:unified,ban-bas:exact,ols:laguerre}, and Hermite type  \cite{mac:hypergeometric,las:hermite,bak-for:calogero-sutherland,die:confluent,xu:orthogonal,dum-ede-shu:MOPS,wad-nis-uji:symmetric,hal-lan:unified}, respectively. Hence, we concentrate on hypergeometric families that (i) are obtained from the Askey-Wilson level via limit transitions rather than parameter specializations and (ii) are endowed with a {\em continuous} orthogonality measure.

For all these families it turns out that the relevant branching polynomials can be conveniently written in terms of expansion coefficients $C^{\mu,n}_{\lambda , r} $
arising from Pieri formulas \cite{die:self-dual,die:properties,die:confluent}
\begin{equation}\label{pieri}
E_r (x_1,\ldots ,x_n) P_\lambda(x_1,\ldots ,x_n)  = \sum_{\substack{\mu\in\Lambda_n \\ \mu\sim_r\lambda}}C^{\mu,n}_{\lambda , r}  P_\mu (x_1,\ldots ,x_n) \qquad (r=1,\ldots , n),
\end{equation}
associated with a suitable choice of generators $E_1 (x_1,\ldots ,x_n) ,\ldots, E_n (x_1,\ldots ,x_n) $ for our algebra of (trigonometric) symmetric
polynomials. The nonvanishing expansion coefficents on the RHS of Eq. \eqref{pieri} are in these cases governed 
by the following proximity relation within $\Lambda_n$:
$\mu\sim_r\lambda$ 
iff there exists a partition $\nu\in\Lambda_n$ with $\nu\subset \lambda$ and $\nu\subset \mu$ such that the skew diagrams $ \lambda /\nu $ and $ \mu / \nu$ are vertical strips with
 $| \lambda /\nu| +| \mu/\nu |\leq r$. Here $|\cdot |$ denotes the number of boxes of the diagram and (recall) the skew diagram $\lambda/\nu$ is a vertical strip iff $\nu_j\leq\lambda_j\leq\nu_j+1$ ($j=1,\ldots ,n$). 

After recalling---in Section \ref{AW:sec}---the explicit branching polynomials from \cite{die-ems:branching} that enable the recursive construction of the Macdonald-Koornwinder polynomials starting from the one-variable Askey-Wilson polynomials, we will work our way down Askey's scheme and provide the  corresponding branching polynomials for
the Wilson level (Section \ref{W:sec}), the continuous Hahn level (Section \ref{cH:sec}), the Jacobi level (Section \ref{J:sec}), the Laguerre level (Section \ref{L:sec}), and the Hermite level (Section \ref{H:sec}). At the bottom level of the symmetric Hermite polynomials, to date only special instances of the corresponding Pieri coefficients are available in closed form in the literature. Our proof of the branching formula mimics in this situation the proof from  \cite{die-ems:branching} for the Macdonald-Koornwinder case and relies on the Hermite degeneration of a (dual) Cauchy identity due to Mimachi detailed in the  appendix at the end of the paper.

{\em Notation.} {\em i)} For future reference, we associated to any pair of partitions $\lambda,\mu\in\Lambda_n$  the following subsets of $\{ 1,\ldots ,n\}$:
$$
J=J(\lambda,\mu):=\{ 1\leq j\leq n\mid \lambda_j\neq \mu_j\},\quad J^c=J^c(\lambda,\mu):=\{ 1\leq j\leq n\mid \lambda_j=\mu_j\} ,$$
$$
 J_+=J_+(\lambda,\mu) :=\{ 1\leq j\leq n\mid \mu_j>\lambda_j\} ,\quad  J_-=J_-(\lambda ,\mu):=\{ 1\leq j\leq n\mid \mu_j<\lambda_j\} ,
$$
and we also define
$$
\epsilon_j=\epsilon_j(J_+,J_-):=\begin{cases} 1 &\text{if} \ j\in J_+, \\ -1&\text{if}\ j\in J_-, \\ 0& \text{otherwise} .\end{cases}
$$
It is immediate from these definitions that if $\mu\sim_r\lambda$, then the cardinality $|J|$ of $J=J(\lambda,\mu)$ is at most $r$ and $\epsilon_j=\mu_j-\lambda_j$ for $j=1,\ldots,n$.

{\em ii)} Following standard conventions, shifted factorials and their $q$-versions are denoted by:
$$(a)_k:=(a)(a+1)\cdots (a+k-1) \quad\text{and}\quad 
(a;q)_k:=(1-a)(1-aq)\cdots (1-aq^{k-1}) ,$$
with $(a)_0=(a;q)_0:= 1$, and
$$(a_1,\ldots,a_l)_k:=(a_1)_k\cdots (a_l)_k,\quad
(a_1,\ldots,a_l;q)_k:=(a_1;q)_k\cdots (a_l;q)_k.$$

{\em iii)} Finally, we will employ the (principal specialization) vectors $\boldsymbol{\tau}=(\tau_1,\ldots ,\tau_n)$ and
$\boldsymbol{\rho}=(\rho_1,\ldots ,\rho_n)$ with components given by
$$\tau_j=t^{n-j}t_0, \qquad \hat{\tau}_j=t^{n-j}\hat{t}_0 $$
and
$$\rho_j=({n-j})g+g_0, \qquad \hat{\rho}_j=({n-j})g+\hat{g}_0 $$
 $(j=1,\ldots ,n)$, where $t$, $t_0$, $\hat{t}_0$ and $g$, $g_0$, $\hat{g}_0$ denote parameters to be specified below.

\section{Askey-Wilson level}\label{AW:sec}

\subsection{Symmetric Macdonald-Koornwinder polynomials \cite{koo:askey-wilson,mac:orthogonal}}
The symmetric Macdonald-Koornwinder polynomials are trigonometric polynomials
\begin{equation}
P_\lambda(x_1,\ldots ,x_n)=P_\lambda(x_1,\ldots ,x_n;q,t,t_l)\qquad ( \lambda\in \Lambda_n)
\end{equation}
determined by the properties in Eqs.
\eqref{pa}, \eqref{pb}, with 
$$M_\lambda(x_1,\ldots ,x_n)= m_\lambda (e^{ix_1}+e^{-ix_1},\ldots ,e^{ix_n}+e^{-ix_n}) ,$$
$$
m_\lambda (z_1,\ldots ,z_n):=\sum_{\nu\in S_n(\lambda)} z_1^{\nu_1}\cdots z_n^{\nu_n} 
$$
(where the sum is over the orbit of $\lambda$ with respect to the action of the symmetric group $S_n$), and 
\begin{align}
\Delta(x_1,\ldots ,x_n) &=\Delta (x_1,\ldots ,x_n;q,t,t_l)  \\
&:=\prod_{1\leq j\leq n} \left| \frac{(e^{2ix_j};q)_\infty }{\prod_{0\leq l\leq 3} (t_l e^{ix_j};q)_\infty}\right|^2
\prod_{1\leq j<k\leq n} \left| \frac{(e^{i(x_j+x_k)},e^{i(x_j-x_k)};q)_\infty }{(te^{i(x_j+x_k)},te^{i(x_j-x_k)};q)_\infty } \right|^2  \nonumber
\end{align}
supported on $\mathcal{D}=[-\pi,\pi]^n$.
Here  it is
assumed that the parameters belong to the domain $0<q, |t|, |t_l| <1$,  with $t$ being real and possibly non-real parameters $t_l$ ($l=0,1,2 ,3$) occurring in complex conjugate pairs.

\subsection{Pieri coefficients \cite{die:self-dual,die:properties,sah:nonsymmetric}}\label{secAW:pieri}
Let $\boldsymbol{e}_m(z_1,\ldots,z_n)$ and $\boldsymbol{h}_m(z_1,\ldots,z_n)$ denote the elementary and the complete symmetric polynomials of degree $m$:
$$
\boldsymbol{e}_m(z_1,\ldots,z_n):=\sum_{1\leq j_1<\cdots <j_m\leq n}
z_{j_1}\cdots z_{j_m} ,
$$
$$
\boldsymbol{h}_m(z_1,\ldots,z_n):=\sum_{1\leq j_1\leq\cdots \leq j_m\leq n}
z_{j_1}\cdots z_{j_m},
$$
with the convention that $\boldsymbol{e}_0=\boldsymbol{h}_0\equiv 1$, and let
\begin{equation*}
\hat{t}_0^2=q^{-1}t_0t_1t_2t_3,\qquad \hat{t}_0\hat{t}_l=t_0t_l\quad (l=1,2,3) .
\end{equation*}

The Pieri coefficients $C^{\mu, n}_{\lambda , r}$ \eqref{pieri} for the Macdonald-Koornwinder polynomials
 associated with the multiplication of $P_\lambda(x_1,\ldots ,x_n;q,t,t_l)$ by
\begin{align}\label{column-ip}
&E_r(x_1,\ldots ,x_n)=E_r(x_1,\ldots ,x_n;   t,t_0)\\
&:=\sum_{0\leq m\leq r} (-1)^{r+m} \boldsymbol{e}_m(e^{ix_1}+e^{-ix_1},\ldots ,e^{ix_n}+e^{-ix_n}) \boldsymbol{h}_{r-m} (\tau_r+\tau_r^{-1},\ldots, \tau_n+\tau_n^{-1})
\nonumber
\end{align}
are given by 
\begin{equation}
C^{\mu, n}_{\lambda , r}=C^{\mu,n}_{\lambda , r}(q,t,t_l):= 
\frac{p_\lambda(q,t,t_l)}{p_\mu (q,t,t_l)} V^{n}_{J_+,J_-} (\lambda ;q,t,t_l) U^{n}_{J^c,r-|J|} (\lambda ; q,t,t_l) ,
\end{equation}
where
\begin{eqnarray*}
\lefteqn{p_\lambda(q,t,t_l):=}&& \\
&& \prod_{1\leq j\leq n}  \frac{\prod_{0\leq l\leq 3} (\hat{t}_l\hat{\tau}_j ;q)_{\lambda_j }}{\tau_j^{\lambda_j} (\hat{\tau}_j^2 ;q)_{2\lambda_j}} 
\prod_{1\leq j<k\leq n} \frac{(t\hat{\tau}_j\hat{\tau}_k;q)_{\lambda_j+\lambda_k} (t\hat{\tau}_j\hat{\tau}_k^{-1};q)_{\lambda_j-\lambda_k}}{(\hat{\tau}_j\hat{\tau}_k;q)_{\lambda_j+\lambda_k}(\hat{\tau}_j\hat{\tau}_k^{-1};q)_{\lambda_j-\lambda_k}}  
\end{eqnarray*}
(which corresponds to the principal specialization value of $P_\lambda(x_1,\ldots ,x_n;q,t,t_l)$),
\begin{align*}
V^{n}_{J_+,J_-} (\lambda ;q,t,t_l)&:=\prod_{j\in J} \frac{\prod_{0\leq l\leq 3}   (1-\hat{t}_l\hat{\tau}_j^{\epsilon_j}q^{\epsilon_j\lambda_j})}{t_0(1-\hat{\tau}_j^{2\epsilon_j}q^{2\epsilon_j\lambda_j})(1-\hat{\tau}_j^{2\epsilon_j}q^{2\epsilon_j\lambda_j+1})} \\
&\times \prod_{\substack{j,j^\prime\in J\\ j<j^\prime}}
\frac{(1-t\hat{\tau}_j^{\epsilon_j}\hat{\tau}_{j^\prime}^{\epsilon_{j^\prime}}q^{\epsilon_j\lambda_j+\epsilon_{j^\prime}\lambda_{j^\prime}})(1-t\hat{\tau}_j^{\epsilon_j}\hat{\tau}_{j^\prime}^{\epsilon_{j^\prime}}q^{\epsilon_j\lambda_j+\epsilon_{j^\prime}\lambda_{j^\prime}+1})}{t(1-\hat{\tau}_j^{\epsilon_j}\hat{\tau}_{j^\prime}^{\epsilon_{j^\prime}}q^{\epsilon_j\lambda_j+\epsilon_{j^\prime}\lambda_{j^\prime}})(1-\hat{\tau}_j^{\epsilon_j}\hat{\tau}_{j^\prime}^{\epsilon_{j^\prime}}q^{\epsilon_j\lambda_j+\epsilon_{j^\prime}\lambda_{j^\prime}+1})} \\
&\times \prod_{j\in J,k\in J^c} \frac{(1-t\hat{\tau}_j^{\epsilon_j}\hat{\tau}_k q^{\epsilon_j\lambda_j+\lambda_k})(1-t\hat{\tau}_j^{\epsilon_j}\hat{\tau}_k^{-1}q^{\epsilon_j\lambda_j-\lambda_k})}{t(1-\hat{\tau}_j^{\epsilon_j}\hat{\tau}_k q^{\epsilon_j\lambda_j+\lambda_k})(1-\hat{\tau}_j^{\epsilon_j}\hat{\tau}_k^{-1}q^{\epsilon_j\lambda_j-\lambda_k})}
\end{align*}
with  $\epsilon_j\equiv \epsilon_j(J_+,J_-)$, and
\begin{align*}
 U^{n}_{K,p} (\lambda ;q,t,t_l)&:=  (-1)^p \sum_{\substack{I_+,I_-\subset K\\ I_+\cap I_-=\emptyset \\ |I_+|+|I_-|= p }} 
 \Biggl( \prod_{j\in I} \frac{\prod_{0\leq l\leq 3}   (1-\hat{t}_l\hat{\tau}_j^{\epsilon_j}q^{\epsilon_j\lambda_j})}{t_0(1-\hat{\tau}_j^{2\epsilon_j}q^{2\epsilon_j\lambda_j})(1-\hat{\tau}_j^{2\epsilon_j}q^{2\epsilon_j\lambda_j+1})} \\
&\times \prod_{\substack{j,j^\prime\in I\\ j<j^\prime}}
\frac{(1-t\hat{\tau}_j^{\epsilon_j}\hat{\tau}_{j^\prime}^{\epsilon_{j^\prime}}q^{\epsilon_j\lambda_j+\epsilon_{j^\prime}\lambda_{j^\prime}})(1-t^{-1}\hat{\tau}_j^{\epsilon_j}\hat{\tau}_{j^\prime}^{\epsilon_{j^\prime}}q^{\epsilon_j\lambda_j+\epsilon_{j^\prime}\lambda_{j^\prime}+1})}{(1-\hat{\tau}_j^{\epsilon_j}\hat{\tau}_{j^\prime}^{\epsilon_{j^\prime}}q^{\epsilon_j\lambda_j+\epsilon_{j^\prime}\lambda_{j^\prime}})(1-\hat{\tau}_j^{\epsilon_j}\hat{\tau}_{j^\prime}^{\epsilon_{j^\prime}}q^{\epsilon_j\lambda_j+\epsilon_{j^\prime}\lambda_{j^\prime}+1})} \\
&\times \prod_{j\in I,k\in K\setminus I} \frac{(1-t\hat{\tau}_j^{\epsilon_j}\hat{\tau}_k q^{\epsilon_j\lambda_j+\lambda_k})(1-t\hat{\tau}_j^{\epsilon_j}\hat{\tau}_k^{-1}q^{\epsilon_j\lambda_j-\lambda_k})}{t(1-\hat{\tau}_j^{\epsilon_j}\hat{\tau}_k q^{\epsilon_j\lambda_j+\lambda_k})(1-\hat{\tau}_j^{\epsilon_j}\hat{\tau}_k^{-1}q^{\epsilon_j\lambda_j-\lambda_k})} \Biggr)
\end{align*}
with  $\epsilon_j\equiv \epsilon_j(I_+,I_-)$, $I=I_+\cup I_-$ and $p=0,\ldots ,|K|$.
In these formulas  it is assumed---by convention---that
 $V^{n}_{J_+,J_-} (\lambda ;q,t,t_l)=1$ if $J$ is empty and that 
 $U^{n}_{K,p} (\lambda ;q,t,t_l)=1$ if $p=0$; furthermore, the index sets  $J=J(\lambda,\mu)$, $J^c=J^c(\lambda,\mu)$ and $J_\pm =J_\pm (\lambda,\mu)$, the signs $\epsilon_j=\epsilon_j(J_+,J_-)$, and the principal specialization vectors $\tau_j$, $\hat{\tau}_j$ ($j=1,\ldots ,n$), are all defined in accordance with the conventions detailed at the end of the introduction.

\begin{remark}
It was observed in \cite[Sec. 5]{kom-nou-shi:kernel} that the polynomial $E_r(x_1,\ldots ,x_n;   t,t_0)$
\eqref{column-ip} amounts to a special instance of Okounkov's hyperoctahedral interpolation polynomial
\cite{oko:bc-type,rai:bcn-symmetric} corresponding to the partition that consists of a single column of size $r$.
\end{remark}

\subsection{Branching polynomials \cite{die-ems:branching}}
For $\lambda\in \Lambda_{n+1}$ and $\mu\in \Lambda_n$ with $\mu\preceq\lambda$, the branching
polynomial $P_{\lambda/\mu}(x)=P_{\lambda/\mu}(x; q ,t ,t_l)$---which arises as the expansion coefficient for the
Macdonald-Koornwinder polynomial in $(n+1)$ variables  in terms of the $n$-variable polynomials  \eqref{bf}---is given explicitly by
\begin{subequations}
\begin{equation}\label{br-pol}
\boxed{P_{\lambda/\mu}(x; q ,t ,t_l)=
\sum_{0\leq k\leq d} B_{\lambda /\mu}^k (q,t,t_l)
\langle x; t_0\rangle_{q,k} }
\end{equation}
where
\begin{equation*}
\langle x; t_0\rangle_{q,k}:=
\prod_{1\leq j\leq k}
(e^{ix}+e^{-ix}-q^{j-1}t_0-q^{-(j-1)}t_0^{-1})
 \end{equation*}
(with $\langle x; t_0\rangle_{q,0}:=1$),
$$d=d(\lambda,\mu):=|\{ 1\leq j\leq m \mid \lambda_j^\prime=\mu_j^\prime +1\}| ,$$ and
\begin{equation}\label{br-coef}
\boxed{B_{\lambda/\mu}^k (q,t,t_l)=
 (-1)^{k+|\lambda|-|\mu|} C^{(n+1)^m-\lambda^\prime,m}_{n^m-\mu^\prime,m-k}(t,q,t_l) \quad (k=0,\ldots ,d).}
 \end{equation}
 \end{subequations}
Here  $|\lambda|=\lambda_1+\cdots +\lambda_n,$ $m=\ell(\lambda^\prime)$ ($=\lambda_1$), and $\lambda^\prime$ ($\in \Lambda_m$) denotes the conjugate partition of $\lambda$, while
$m^n-\mu$ with $\mu\subset m^n$ stands for the partition such that $(m^n-\mu)_j=m-\mu_{n+1-j}$  $(j=1,\ldots,n). $

Apart from exploiting the above Pieri formulas,
the proof in \cite{die-ems:branching} of this branching rule is based on
Mimachi's Cauchy identity for the Macdonald-Koornwinder polynomials \cite[Thm. 2.1]{mim:duality} (cf. Eq. \eqref{AW:cauchy} below) as well as on a special `column-row' case 
 \cite[Lem. 5.1]{kom-nou-shi:kernel}
of the  Cauchy identity for Okounkov's hyperoctahedral interpolation polynomials \cite[Thm. 6.2]{oko:bc-type} with shifted variables
in accordance with \cite[Thm. 3.16]{rai:bcn-symmetric} (cf. Remark \ref{br-proof:rem} at the end of Section \ref{H:sec} below).

\begin{remark}\label{mac-br:rem}
It is known that the highest-degree leading homogeneous terms of the Macdonald-Koornwinder polynomial $P_\lambda(x_1,\ldots ,x_n;q,t,t_l)$ consist of the (monic) Macdonald polynomial  $P_\lambda(x_1,\ldots ,x_n;q,t)$  \cite[\S 5.2]{die:commuting}.
This ties in with the recent observation in  \cite[Rem. 6.1]{koo:okounkov} that---by filtering the terms of leading degree on both sides of
the above branching formula for the Macdonald-Koornwinder polynomials---a celebrated branching rule for the Macdonald polynomials 
\cite{mac:symmetric,las-war:branching,sun:representation} is recovered:
$$
P_\lambda (x_1,\ldots ,x_n,x;q,t)=
\sum_{\substack{\mu\in\Lambda_n,\,\mu\subset\lambda \\ \lambda/\mu \ \text{horizontal\ strip}}}  P_\mu (x_1,\ldots ,x_n;q,t )
P_{\lambda/\mu}(x;q,t)
$$
$(\lambda\in\Lambda_{n+1})$, where
$$
P_{\lambda/\mu}(x;q,t)= e^{ix ( |\lambda |-|\mu |)}   B_{\lambda/\mu}(q,t)
$$
and
\begin{align*}
B_{\lambda/\mu}(q,t)=&B^{d}_{\lambda/\mu}(q,t,t_l)=C^{(n+1)^m-\lambda^\prime,m}_{n^m-\mu^\prime ,m-d}(t,q,t_l)
\quad\text{with} \ d=|\lambda|-|\mu|
 \\
= & \frac{p_{n^m-\mu^\prime}(t,q,t_l)}{p_{(n+1)^m-\lambda^\prime}(t,q,t_l)}
V^{m}_{J(n^m-\mu^\prime ,(n+1)^m-\lambda^\prime),\emptyset} (n^m-\mu^\prime ;t,q,t_l )  \\
=&
\prod_{1\leq j<k\leq m} \frac{(q^{1+k-j};t)_{\mu_j^\prime-\mu_k^\prime}}{(q^{k-j};t)_{\mu_j^\prime-\mu_k^\prime}} 
 \frac{(q^{k-j};t)_{\lambda_j^\prime-\lambda_k^\prime}}{(q^{1+k-j};t)_{\lambda_j^\prime-\lambda_k^\prime}} \\
 &\times \prod_{\substack{ 1\leq j,k\leq m\\ \mu_j^\prime\neq \lambda_j^\prime \\ \mu_k^\prime = \lambda_k^\prime}}
 \biggl( \frac{1-q^{1+k-j}t^{\mu_j^\prime-\mu_k^\prime} }{1-q^{k-j}t^{\mu_j^\prime-\mu_k^\prime} }\biggr) 
  \prod_{\substack{ 1\leq j<k\leq m\\ \mu_j^\prime= \lambda_j^\prime \\ \mu_k^\prime \neq \lambda_k^\prime}} q^{-1} 
 \\
 =&
  \prod_{\substack{ 1\leq j<k\leq m\\ \mu_j^\prime= \lambda_j^\prime \\ \mu_k^\prime \neq \lambda_k^\prime}}
   \biggl( \frac{1-q^{1+k-j}t^{\lambda_j^\prime-\lambda_k^\prime} }{1-q^{k-j}t^{\lambda_j^\prime-\lambda_k^\prime} }\biggr) 
 \biggl( \frac{1-q^{-1+k-j}t^{\mu_j^\prime-\mu_k^\prime} }{1-q^{k-j}t^{\mu_j^\prime-\mu_k^\prime} }\biggr) 
\\
\stackrel{*}{=}&
 \prod_{1\leq j\leq k\leq \ell (\mu)} 
 \frac{(q^{\mu_j-\mu_k}t^{1+k-j},q^{1+\mu_j-\lambda_{k+1}} t^{k-j};q)_{\lambda_j-\mu_j}}{(q^{1+\mu_j-\mu_k}t^{k-j},q^{\mu_j-\lambda_{k+1}}t^{1+k-j};q)_{\lambda_j-\mu_j}} 
\end{align*}
* cf. Eq. (6.13), Rem. 2., and Example 2.(b) of \cite[Ch. VI.6]{mac:symmetric}.
\end{remark}

\subsection{Whittaker limit ($t\to 0$)}
For $t\to 0$ the Macdonald-Koornwinder polynomial degenerates into a deformed $q$-Whittaker function $P_\lambda(x_1,\ldots ,x_n;q,0,t_l)$ that diagonalizes Ruijsenaars's $q$-difference Toda chain with one-sided integrable boundary interactions of Askey-Wilson type
\cite{die-ems:integrable}. The corresponding branching polynomial for $P_\lambda(x_1,\ldots ,x_n;q,0,t_l)$---obtained from
Eq. \eqref{br-pol}  in the limit $t\to 0$---is given by
\begin{subequations}
\begin{equation}\label{br-pol-whittaker}
\boxed{P_{\lambda/\mu}(x; q ,0,t_l)=
\sum_{0\leq k\leq d} B_{\lambda /\mu}^k (q,0,t_l)
\langle x; t_0\rangle_{q,k} }
\end{equation}
with coefficients
\begin{equation}
\boxed{B_{\lambda/\mu}^k (q,0,t_l)=
 (-1)^{k+|\lambda|-|\mu|} C^{(n+1)^m-\lambda^\prime,m}_{n^m-\mu^\prime,m-k}(0,q,t_l)}
 \end{equation}
 \end{subequations}
governed by the Pieri coefficients for the Macdonald-Koornwinder polynomials at $q=0$:
\begin{equation}
C^{\mu,n}_{\lambda , r}(0,t,t_l):= 
\frac{p_\lambda(0,t,t_l)}{p_\mu (0,t,t_l)} V^{n}_{J_+,J_-} (\lambda ;0,t,t_l) U^{n}_{J^c,r-|J|} (\lambda ; 0,t,t_l) .
\end{equation}
These Pieri coefficients are given explicitly by \cite{die-ems:quantum}:
\begin{align*}
&p_\lambda(0,t,t_l):= \\
& \prod_{\substack{1\leq j\leq n\\ \lambda_j>0}}  \Bigl( \tau_j^{-\lambda_j}{\prod_{1\leq l\leq 3} (1-t_0t_lt^{n-j})}\Bigr)
\prod_{\substack{1\leq j\leq n\\ \lambda_j=1}}  (1-t_0t_1t_2t_3t^{n-j+m_0(\lambda)})^{-1}
\prod_{\substack{1\leq j<k\leq n\\ \lambda_j>\lambda_k}} \frac{1-t^{1+k-j}}{1-t^{k-j}}  ,
\end{align*}
\begin{align*}
 &V^n_{J_+,J_-}(\lambda ;0,t,t_l) :=\\
 & \prod_{\substack{j\in J_+  \\ \lambda_j=0}} 
\frac{ (1-t_0t_1t_2t_3 t^{n-j+m_0(\lambda )+m_1(\lambda)-m_1^+(\lambda)})
\prod_{1\leq l\leq 3}(1-t_0t_l t^{n-j}) } {  (1-t_0t_1t_2t_3 t^{2(n-j)}) (1-t_0t_1t_2t_3 t^{2(n-j)+1}) } \nonumber \\
& \times \prod_{\substack{j\in J_+  \\ \lambda_j=1}} (1-t_0t_1t_2t_3 t^{n-j+m_0(\lambda)}) \prod_{\substack{j\in J_-  \\ \lambda_j=1}} 
\frac{ (1-t_0t_1t_2t_3 t^{n-j-1})\prod_{1\leq l < m\leq 3}(1-t_lt_mt^{n-j}) }{  (1-t_0t_1t_2t_3 t^{2(n-j)}) (1-t_0t_1t_2t_3 t^{2(n-j)-1}) }\nonumber  \\
&\times \prod_{\substack{1\leq j < k\leq n\\  \lambda_j=\lambda_k, \epsilon_j>\epsilon_k}}
\frac{1-t^{1+k-j}}{1-t^{k-j}}    \prod_{1\leq j\leq n} \tau_j^{-\epsilon_j} , \nonumber
\end{align*} 
with $\epsilon_j\equiv\epsilon_j(J_+,J_-)$, and
 \begin{align*}\label{Ukm}
&U^n_{K, p}(\lambda ;0,t,t_l) := (-1)^p \times \\
&\sum_{\substack{I_+,I_-\subset K\\ I_+\cap I_-=\emptyset \\ |I_+|+|I_-|= p }} 
\Biggl(
\prod_{\substack{j\in I_+\\ \lambda_j=0}}   \frac{\prod_{1\leq l\leq 3}(1-t_0t_l t^{n-j}) }{1-t_0t_1t_2t_3 t^{2(n-j)}}
\prod_{\substack{j\in I_+\\ \lambda_j=1}}  (1-t_0t_1t_2t_3 t^{n-j}) \\
&\times
 \prod_{\substack{j\in I_-\\\lambda_j=1}}\frac{\prod_{1\leq l <m\leq 3} (1- t_lt_mt^{n-j})}{1-t_0t_1t_2t_3 t^{2(n-j)}}
\prod_{\substack{j , k\in K \\  \lambda_j=\lambda_k, \epsilon_j>\epsilon_k}}
\frac{1-t^{1+k-j}}{1-t^{k-j}}
\prod_{\substack{j\in I_-,\, k\in I_+\\  \lambda_j=\lambda_k+1}}
\frac{1-t^{1+k-j}}{1-t^{k-j}}
\nonumber \\
 & \times   \prod_{\substack{j,k\in K, j<k \\ \epsilon_j+\epsilon_k\in \{ -2,1,2\} \\ \lambda_j=1, \lambda_k=\delta_{1+\epsilon_k}}}\frac{1-t_0t_1t_2t_3 t^{2n+1-j-k}}{1-t_0t_1t_2t_3 t^{2n-j-k}}
  \prod_{\substack{j\in I_+\cup I_-,k\in K\setminus I_- \\ j<k, \epsilon_k-\epsilon_j\in\{ 0,1\}\\ \lambda_j=\delta_{1+\epsilon_j}, \lambda_k=0 }} \frac{1-t_0t_1t_2t_3 t^{2n-1-j-k}}{1-t_0t_1t_2t_3 t^{2n-j-k}}
 \nonumber \\
&\times \prod_{j\in K} t_0^{-\epsilon_j}  
\prod_{\substack{j,k\in K,\, j <k \\ \epsilon_j\neq\epsilon_k=0}} t^{-\epsilon_j}
\prod_{\substack{j,k\in K,\, j <k \\ \lambda_j=\lambda_k,\, \epsilon_k-\epsilon_j=1}} t^{-1}
\Biggr) ,
\nonumber
 \end{align*} 
with $\epsilon_j\equiv\epsilon_j(I_+,I_-)$, where we have employed the additional notation
$m_l (\lambda):=|\{ 1\leq j\leq n \mid \lambda_j=l\} |$,  $m_l^+ (\lambda):=|\{ j\in J_+ \mid \lambda_j=l\} |$, and $\delta_m:=1$ if $m=0$ and
$\delta_m:=0$ otherwise.

\begin{remark}
It is evident from these Pieri coefficients that the branching polynomial $P_{\lambda/\mu}(x; q ,0,t_l)$ \eqref{br-pol-whittaker} for the deformed $q$-Whittaker function $P_\lambda(x_1,\ldots ,x_n;q,0,t_l)$ simplifies considerably 
when one or more of the parameters $t_1,t_2$ or $t_3$ vanish. From the perspective of Ruijsenaars' $q$-difference Toda chain, such parameter reductions correspond to degenerations of the interaction at the boundary
\cite[Sec. 7]{die-ems:integrable}.
\end{remark}

\section{Wilson level}\label{W:sec}

\subsection{Symmetric Wilson polynomials \cite{die:multivariable,die:properties,gro:multivariable}}
The symmetric Wilson polynomials are even polynomials
\begin{equation}
P_\lambda(x_1,\ldots ,x_n)=P_\lambda^{W}(x_1,\ldots ,x_n; g,g_l)\qquad ( \lambda\in \Lambda_n)
\end{equation}
of the form in Eqs.
\eqref{pa}, \eqref{pb}, with
$$M_\lambda(x_1,\ldots ,x_n)= m_\lambda (x_1^2,\ldots ,x_n^2) $$
and 
\begin{align}
&\Delta(x_1,\ldots ,x_n) =\Delta^W(x_1,\ldots ,x_n; g,g_l)  \\
&:=\prod_{1\leq j\leq n} \left| \frac{ \prod_{0\leq l\leq 3} \Gamma(g_l+ix_j)  }
{ \Gamma(2i x_j) }
\right|^2
\prod_{1\leq j<k\leq n} \left| \frac{ \Gamma(g+i(x_j+x_k)) \Gamma(g+i(x_j-x_k)) }{\Gamma(i(x_j+ x_k)) \Gamma(i(x_j- x_k))}
\right|^2  \nonumber
\end{align}
supported on $\mathcal{D}=\mathbb{R}^n$.
Here $\Gamma (\cdot )$ denotes the gamma function and it is
assumed that $g, \text{Re}(g_l)>0$, with possibly non-real parameters $g_l$ occurring in complex conjugate pairs ($l=0,1,2 ,3$).

\subsection{Pieri coefficients \cite{die:properties}}
Let
\begin{equation*}
\hat{g}_0=\frac{1}{2}(g_0+g_1+g_2+g_3-1)\quad  \text{and}\quad \hat{g}_0+\hat{g}_l=g_0+g_l\quad (l=1,2,3) .
\end{equation*}
The Pieri coefficients $C^{\mu, n}_{\lambda , r}$ \eqref{pieri} for the symmetric Wilson polynomials
 associated with the multiplication of $P^W_\lambda(x_1,\ldots ,x_n;g,g_l)$ by
\begin{align}
&E_r(x_1,\ldots ,x_n)=E_r^W(x_1,\ldots ,x_n;   g,g_0)\\
&:=(-1)^{r}\sum_{0\leq m\leq r}  \boldsymbol{e}_m(x_1^2,\ldots ,x_n^2) \boldsymbol{h}_{r-m} (\rho_r^2,\ldots, \rho_n^2)
\nonumber
\end{align}
are given by 
\begin{equation}
C^{\mu, n}_{\lambda , r}=C^{W,\mu,n}_{\lambda , r}(g,g_l):= 
\frac{p^W_\lambda(g,g_l)}{p^W_\mu (g,g_l)} V^{W,n}_{J_+,J_-} (\lambda ;g,g_l) U^{W,n}_{J^c,r-|J|} (\lambda ; g,g_l) ,
\end{equation}
where
\begin{align*}
&p^W_\lambda(g,g_l):= \\
&(-1)^{|\lambda|} \prod_{1\leq j\leq n}  \frac{\prod_{0\leq l\leq 3} (\hat{g}_l+\hat{\rho}_j )_{\lambda_j }}{ (2\hat{\rho}_j)_{2\lambda_j}} 
\prod_{1\leq j<k\leq n} \frac{(g+\hat{\rho}_j+\hat{\rho}_k)_{\lambda_j+\lambda_k} (g+\hat{\rho}_j-\hat{\rho}_k)_{\lambda_j-\lambda_k}}{(\hat{\rho}_j+\hat{\rho}_k)_{\lambda_j+\lambda_k}(\hat{\rho}_j-\hat{\rho}_k)_{\lambda_j-\lambda_k}}  ,
\end{align*}
\begin{align*}
&V^{W,n}_{J_+,J_-} (\lambda ;g,g_l):=\prod_{j\in J} \frac{\prod_{0\leq l\leq 3}   (\hat{g}_l+\epsilon_j(\hat{\rho}_j+\lambda_j))}{2\epsilon_j(\hat{\rho}_j+\lambda_j)(1+2\epsilon_j(\hat{\rho}_j+\lambda_j))} \\
&\times \prod_{\substack{j,j^\prime\in J\\ j<j^\prime}}
\frac{(g+\epsilon_j(\hat{\rho}_j+\lambda_j)+\epsilon_{j^\prime}(\hat{\rho}_{j^\prime
}+\lambda_{j^\prime}))(1+g+\epsilon_j(\hat{\rho}_j+\lambda_j)+\epsilon_{j^\prime}(\hat{\rho}_{j^\prime
}+\lambda_{j^\prime}))}{(\epsilon_j(\hat{\rho}_j+\lambda_j)+\epsilon_{j^\prime}(\hat{\rho}_{j^\prime
}+\lambda_{j^\prime}))(1+\epsilon_j(\hat{\rho}_j+\lambda_j)+\epsilon_{j^\prime}(\hat{\rho}_{j^\prime
}+\lambda_{j^\prime}))} \\
&\times \prod_{j\in J,k\in J^c} \frac{(g+\epsilon_j(\hat{\rho}_j+\lambda_j)+\hat{\rho}_k +\lambda_k)(g+\epsilon_j(\hat{\rho}_j+\lambda_j)-\hat{\rho}_k -\lambda_k)}{(\epsilon_j(\hat{\rho}_j+\lambda_j)+\hat{\rho}_k +\lambda_k)(\epsilon_j(\hat{\rho}_j+\lambda_j)-\hat{\rho}_k -\lambda_k)}
\end{align*}
with $\epsilon_j\equiv \epsilon_j(J_+,J_-)$, and
\begin{align*}
& U^{W,n}_{K,p} (\lambda ;g,g_l):=  (-1)^p \sum_{\substack{I_+,I_-\subset K\\ I_+\cap I_-=\emptyset \\  |I_+|+|I_-|= p }} 
 \Biggl( \prod_{j\in I}  \frac{\prod_{0\leq l\leq 3}   (\hat{g}_l+\epsilon_j(\hat{\rho}_j+\lambda_j))}{2\epsilon_j(\hat{\rho}_j+\lambda_j)(1+2\epsilon_j(\hat{\rho}_j+\lambda_j))} \\
&\times \prod_{\substack{j,j^\prime\in I\\ j<j^\prime}}
\frac{(g+\epsilon_j(\hat{\rho}_j+\lambda_j)+\epsilon_{j^\prime}(\hat{\rho}_{j^\prime
}+\lambda_{j^\prime}))(1-g+\epsilon_j(\hat{\rho}_j+\lambda_j)+\epsilon_{j^\prime}(\hat{\rho}_{j^\prime
}+\lambda_{j^\prime}))}{(\epsilon_j(\hat{\rho}_j+\lambda_j)+\epsilon_{j^\prime}(\hat{\rho}_{j^\prime
}+\lambda_{j^\prime}))(1+\epsilon_j(\hat{\rho}_j+\lambda_j)+\epsilon_{j^\prime}(\hat{\rho}_{j^\prime
}+\lambda_{j^\prime}))}  \\
&\times \prod_{j\in I,k\in K\setminus I} \frac{(g+\epsilon_j(\hat{\rho}_j+\lambda_j)+\hat{\rho}_k +\lambda_k)(g+\epsilon_j(\hat{\rho}_j+\lambda_j)-\hat{\rho}_k -\lambda_k)}{(\epsilon_j(\hat{\rho}_j+\lambda_j)+\hat{\rho}_k +\lambda_k)(\epsilon_j(\hat{\rho}_j+\lambda_j)-\hat{\rho}_k -\lambda_k)}\Biggr)
\end{align*}
with $\epsilon_j\equiv \epsilon_j(I_+,I_-)$, $I=I_+\cup I_-$ and  $p=0,\ldots ,|K|$. Here $\rho_j$ and $\hat{\rho}_j$ ($j=1,\ldots ,n$) are defined in accordance with the conventions at the end of the introduction.

\subsection{Branching polynomials}
Upon performing a rescaling of the trigonometric variables $x_j\to \alpha x_j$ and picking parameters of the form
\begin{subequations}
\begin{equation}\label{AW->Wa}
q=e^{-\alpha},\quad t=e^{-\alpha g}, \quad  t_l=e^{-\alpha g_l}\ \ (l=0,1,2,3),
\end{equation}
the Macdonald-Koornwinder polynomials degenerate into the symmetric Wilson polynomials in the rational limit $\alpha\to 0$
\cite[Sec. 4.1]{die:properties}:
\begin{equation}\label{AW->Wb}
P_\lambda^{W}(x_1,\ldots ,x_n; g,g_l)=\lim_{\alpha\to 0} 
(-\alpha^{-2})^{|\lambda|}
P_\lambda (\alpha x_1,\ldots ,\alpha x_n; q,t,t_l) .
\end{equation}
\end{subequations}

The corresponding branching polynomials $(-\alpha^{2})^{|\lambda |-|\mu |}P_{\lambda/\mu}(\alpha x; q,t,t_l)$ 
\eqref{br-pol}, \eqref{br-coef} are seen to
degenerate in this limit into the following branching polynomials for the symmetric Wilson polynomials:
\begin{subequations}
\begin{equation}\label{br-pol-W}
\boxed{P^W_{\lambda/\mu}(x; g,g_l) =
\sum_{0\leq k\leq d} B_{\lambda/\mu}^{W,k} (g,g_l)
(g_0+ix,g_0-ix)_k}
\end{equation}
with
\begin{equation}\label{br-coef-W}
\boxed{B_{\lambda/\mu}^{W,k} (g,g_l)=
(-1)^{|\lambda|-|\mu|+m}g^{2 (|\lambda|-|\mu|-k)} 
C^{W,(n+1)^m-\lambda^\prime,m}_{n^m-\mu^\prime,m-k}({\textstyle \frac{1}{g},\frac{g_l}{g}})}
 \end{equation}
 \end{subequations}
for $k=0,\ldots ,d=d(\lambda,\mu)$. Here one uses that
$$
\lim_{\alpha \to 0} (-\alpha^{-2})^k \langle \alpha x, t_0\rangle_{q,k}=(g_0+ix)_k (g_0-ix)_k
$$
and that
$$
\lim_{\alpha\to 0} (-\alpha^{-2})^{|\lambda|- |\mu|} \alpha^{-2r} C^{\mu,n}_{\lambda,r}(q,t,t_l)
=C^{W,\mu,n}_{\lambda,r}(g,g_l) .
$$

\section{Continuous Hahn level}\label{cH:sec}

\subsection{Symmetric continuous Hahn polynomials
\cite{die:multivariable,die:properties}}

The symmetric continuous Hahn polynomials
\begin{equation}
P_\lambda(x_1,\ldots ,x_n)=P_\lambda^{cH}(x_1,\ldots ,x_n; g,g_0,g_1)\qquad ( \lambda\in \Lambda_n)
\end{equation}
are of the form in Eqs.
\eqref{pa}, \eqref{pb}, with
$$M_\lambda(x_1,\ldots ,x_n)= m_\lambda (x_1,\ldots ,x_n) $$
and 
\begin{align}
&\Delta(x_1,\ldots ,x_n) =\Delta^{cH}(x_1,\ldots ,x_n; g,g_0,g_1) \\
&:= \prod_{1\leq j\leq n} \left| \Gamma(g_0+ix_j) \Gamma(g_1+ix_j) 
\right|^2
\prod_{1\leq j<k\leq n} \left| \frac{ \Gamma(g+i(x_j-x_k)) }{\Gamma(i(x_j- x_k))}
\right|^2  \nonumber
\end{align}
supported on $\mathcal{D}=\mathbb{R}^n$. Here it is
assumed that $g$ and $ \text{Re}(g_0)$, $ \text{Re}(g_1)$ are all positive.

\subsection{Pieri coefficients \cite{die:properties}}
Let
$$
\hat{g}_0=-\frac{1}{2}+\text{Re}(g_0)+\text{Re}(g_1),\qquad \hat{g}_1=\frac{1}{2}+\text{Re}(g_0)-\text{Re}(g_1)
$$
and
$$
\hat{g}_2=\frac{1}{2}+i(\text{Im}(g_0)-\text{Im}(g_1)).
$$
The Pieri coefficients $C^{\mu, n}_{\lambda , r}$ \eqref{pieri} for the symmetric continuous Hahn polynomials
 associated with the multiplication of $P^{cH}_\lambda(x_1,\ldots ,x_n;g,g_0,g_1)$ by
\begin{align}
&E_r(x_1,\ldots ,x_n)=E_r^{cH}(x_1,\ldots ,x_n;   g,g_0)\\
&:=(-1)^{r} \sum_{0\leq m\leq r}  \boldsymbol{e}_m(ix_1,\ldots ,ix_n) \boldsymbol{h}_{r-m} (\rho_r,\ldots, \rho_n)
\nonumber
\end{align}
are given by 
\begin{equation}
C^{\mu, n}_{\lambda , r}=C^{cH,\mu,n}_{\lambda , r}(g,g_0,g_1):= 
\frac{p^{cH}_\lambda(g,g_0,g_1)}{p^{cH}_\mu (g,g_0,g_1)} V^{cH,n}_{J_+,J_-} (\lambda ;g,g_0,g_1) U^{cH,n}_{J^c,r-|J|} (\lambda ; g,g_0,g_1) ,
\end{equation}
where
\begin{align*}
&p^{cH}_\lambda(g,g_0,g_1):= \\
&i^{|\lambda|} \prod_{1\leq j\leq n}  \frac{\prod_{0\leq l\leq 2} (\hat{g}_l+\hat{\rho}_j )_{\lambda_j }}{ (2\hat{\rho}_j)_{2\lambda_j}} 
\prod_{1\leq j<k\leq n} \frac{(g+\hat{\rho}_j+\hat{\rho}_k)_{\lambda_j+\lambda_k} (g+\hat{\rho}_j-\hat{\rho}_k)_{\lambda_j-\lambda_k}}{(\hat{\rho}_j+\hat{\rho}_k)_{\lambda_j+\lambda_k}(\hat{\rho}_j-\hat{\rho}_k)_{\lambda_j-\lambda_k}}  ,
\end{align*}
\begin{align*}
&V^{cH,n}_{J_+,J_-} (\lambda ;g,g_0,g_1):=\prod_{j\in J} \frac{\prod_{0\leq l\leq 2}   (\hat{g}_l+\epsilon_j(\hat{\rho}_j+\lambda_j))}{2\epsilon_j(\hat{\rho}_j+\lambda_j)(1+2\epsilon_j(\hat{\rho}_j+\lambda_j))} \\
&\times \prod_{\substack{j,j^\prime\in J\\ j<j^\prime}}
\frac{(g+\epsilon_j(\hat{\rho}_j+\lambda_j)+\epsilon_{j^\prime}(\hat{\rho}_{j^\prime
}+\lambda_{j^\prime}))(1+g+\epsilon_j(\hat{\rho}_j+\lambda_j)+\epsilon_{j^\prime}(\hat{\rho}_{j^\prime
}+\lambda_{j^\prime}))}{(\epsilon_j(\hat{\rho}_j+\lambda_j)+\epsilon_{j^\prime}(\hat{\rho}_{j^\prime
}+\lambda_{j^\prime}))(1+\epsilon_j(\hat{\rho}_j+\lambda_j)+\epsilon_{j^\prime}(\hat{\rho}_{j^\prime
}+\lambda_{j^\prime}))} \\
&\times \prod_{j\in J,k\in J^c} \frac{(g+\epsilon_j(\hat{\rho}_j+\lambda_j)+\hat{\rho}_k +\lambda_k)(g+\epsilon_j(\hat{\rho}_j+\lambda_j)-\hat{\rho}_k -\lambda_k)}{(\epsilon_j(\hat{\rho}_j+\lambda_j)+\hat{\rho}_k +\lambda_k)(\epsilon_j(\hat{\rho}_j+\lambda_j)-\hat{\rho}_k -\lambda_k)}
\end{align*}
with $\epsilon_j\equiv \epsilon_j(I_+,I_-)$, and
\begin{align*}
 &U^{cH,n}_{K,p} (\lambda ;g,g_0,g_1):=  (-1)^p \sum_{\substack{I_+,I_-\subset K\\ I_+\cap I_-=\emptyset \\ |I_+|+|I_-|= p }} 
 \Biggl( \prod_{j\in I}  \frac{\prod_{0\leq l\leq 2}   (\hat{g}_l+\epsilon_j(\hat{\rho}_j+\lambda_j))}{2\epsilon_j(\hat{\rho}_j+\lambda_j)(1+2\epsilon_j(\hat{\rho}_j+\lambda_j))} \\
&\times \prod_{\substack{j,j^\prime\in I\\ j<j^\prime}}
\frac{(g+\epsilon_j(\hat{\rho}_j+\lambda_j)+\epsilon_{j^\prime}(\hat{\rho}_{j^\prime
}+\lambda_{j^\prime}))(1-g+\epsilon_j(\hat{\rho}_j+\lambda_j)+\epsilon_{j^\prime}(\hat{\rho}_{j^\prime
}+\lambda_{j^\prime}))}{(\epsilon_j(\hat{\rho}_j+\lambda_j)+\epsilon_{j^\prime}(\hat{\rho}_{j^\prime
}+\lambda_{j^\prime}))(1+\epsilon_j(\hat{\rho}_j+\lambda_j)+\epsilon_{j^\prime}(\hat{\rho}_{j^\prime
}+\lambda_{j^\prime}))}  \\
&\times \prod_{j\in I,k\in K\setminus I} \frac{(g+\epsilon_j(\hat{\rho}_j+\lambda_j)+\hat{\rho}_k +\lambda_k)(g+\epsilon_j(\hat{\rho}_j+\lambda_j)-\hat{\rho}_k -\lambda_k)}{(\epsilon_j(\hat{\rho}_j+\lambda_j)+\hat{\rho}_k +\lambda_k)(\epsilon_j(\hat{\rho}_j+\lambda_j)-\hat{\rho}_k -\lambda_k)}\Biggr)
\end{align*}
with $\epsilon_j\equiv \epsilon_j(I_+,I_-)$, $I=I_+\cup I_-$ and  $p=0,\ldots ,|K|$.

\subsection{Branching polynomials}
Upon picking parameters of the form
\begin{subequations}
\begin{equation}\label{AW->cHa}
q=e^{-\alpha},\quad t=e^{-\alpha g}, \quad  t_l=\bar{t}_{l+2}=-ie^{-\alpha g_l}\
(l=0,1),
\end{equation}
the Macdonald-Koornwinder polynomials degenerate into the symmetric continuous Hahn polynomials in the shifted rational limit
\cite[Sec. 4.2]{die:properties}:
\begin{equation}\label{AW->cHb}
P_\lambda^{cH}(x_1,\ldots ,x_n; g,g_0,g_1)=\lim_{\alpha\to 0} 
(2\alpha)^{-|\lambda|}
P_\lambda (\alpha x_1-\frac{\pi}{2},\ldots ,\alpha x_n-\frac{\pi}{2}; q,t,t_l) .
\end{equation}
\end{subequations}

The corresponding branching polynomials $(2\alpha)^{-|\lambda |+|\mu |}P_{\lambda/\mu}(\alpha x-\frac{\pi}{2}; q,t,t_l)$ 
\eqref{br-pol}, \eqref{br-coef} are seen to
degenerate in this limit into the following branching polynomials for the symmetric continuous Hahn polynomials:
\begin{subequations}
\begin{equation}\label{br-pol-cH}
\boxed{P^{cH}_{\lambda/\mu}(x; g,g_0,g_1) =
\sum_{0\leq k\leq d} B_{\lambda/\mu}^{cH,k} (g,g_0, g_1)
(g_0+ix)_k }
\end{equation}
with
\begin{align}\label{br-coef-cH}
\boxed{B_{\lambda/\mu}^{cH,k} (g,g_0,g_1)=
i^m (-1)^{|\lambda|-|\mu|}g^{|\lambda|-|\mu|-k} 
C^{cH,(n+1)^m-\lambda^\prime,m}_{n^m-\mu^\prime,m-k}({\textstyle \frac{1}{g},\frac{g_0}{g}, \frac{g_1}{g} }) }\nonumber\\
 \end{align}
 \end{subequations}
for $k=0,\ldots ,d=d(\lambda,\mu)$. Here one uses that
$$
\lim_{\alpha \to 0} (-2i\alpha)^{-k} \langle \alpha x-\frac{\pi}{2}, t_0\rangle_{q,k}=(g_0+ix)_k$$
and that
$$
\lim_{\alpha\to 0} (2\alpha)^{-|\lambda|+ |\mu|} (2\alpha i)^{-r} C^{\mu,n}_{\lambda,r}(q,t,t_l)
=C^{cH,\mu,n}_{\lambda,r}(g,g_0,g_1) .
$$

\section{Jacobi level}\label{J:sec}

\subsection{Symmetric Jacobi polynomials
\cite{vre:formulas,hec-opd:hypergeometric,hec:hypergeometric,mac:hypergeometric,deb:systeme,las:jacobi,bee-opd:certain,bak-for:calogero-sutherland,die-lap-mor:determinantal,dum-ede-shu:MOPS,hal-lan:unified,koo:okounkov}}

The symmetric Jacobi polynomials are trigonometric polynomials
\begin{equation}
P_\lambda(x_1,\ldots ,x_n)=P_\lambda^{J}(x_1,\ldots ,x_n; g,g_0,g_1)\qquad ( \lambda\in \Lambda_n)
\end{equation}
of the form in Eqs.
\eqref{pa}, \eqref{pb}, with
$$M_\lambda(x_1,\ldots ,x_n)= m_\lambda (e^{ix_1}+e^{-ix_1},\ldots ,e^{ix_n}+e^{-ix_n}) $$
and 
\begin{align}
&\Delta(x_1,\ldots ,x_n) =\Delta^{J}(x_1,\ldots ,x_n; g,g_0,g_1) \\
&:= \prod_{1\leq j\leq n} \left| \sin \Bigl(\frac{x_j}{2}\Bigr) \right|^{2g_0-1} \left| \cos \Bigl(\frac{x_j}{2}\Bigr) \right|^{2g_1-1}
\prod_{1\leq j<k\leq n} \left|  \sin \Bigl(\frac{x_j+x_k}{2}\Bigr)  \sin \Bigl(\frac{x_j-x_k}{2}\Bigr) 
\right|^{2 g} \nonumber
\end{align}
supported on $\mathcal{D}=[-\pi,\pi]^n$. Here it is
assumed that the parameters $g,g_0$ and $g_1$ are all positive.

\subsection{Pieri coefficients \cite{die:properties}}
Let
$$
\hat{g}_0= \frac{1}{2}(g_0+g_1-1)\quad\text{and}\quad \hat{g}_1= \frac{1}{2}(g_0-g_1+1).
$$
The Pieri coefficients $C^{\mu, n}_{\lambda , r}$ \eqref{pieri} for the symmetric Jacobi polynomials
 associated with the multiplication of $P^{J}_\lambda(x_1,\ldots ,x_n;g,g_0,g_1)$ by
\begin{equation}
E_r(x_1,\ldots ,x_n)=E_r^{J}(x_1,\ldots ,x_n)\\
:=(-1)^{r}  \boldsymbol{e}_r\left(\sin^2\left(\frac{x_1}{2}\right) ,\ldots , \sin^2\left(\frac{x_n}{2}\right) \right) 
\end{equation}
are given by 
\begin{equation}
C^{\mu, n}_{\lambda , r}=C^{J,\mu,n}_{\lambda , r}(g,g_0,g_1):= 
\frac{p^{J}_\lambda(g,g_0,g_1)}{p^{J}_\mu (g,g_0,g_1)} V^{J,n}_{J_+,J_-} (\lambda ;g,g_0,g_1) U^{J,n}_{J^c,r-|J|} (\lambda ; g,g_0,g_1) ,
\end{equation}
where
\begin{align*}
&p^{J}_\lambda(g,g_0,g_1):= \\
&4^{|\lambda|} \prod_{1\leq j\leq n}  \frac{(\hat{g}_0+\hat{\rho}_j, \hat{g}_1+\hat{\rho}_j )_{\lambda_j }}{ (2\hat{\rho}_j)_{2\lambda_j}} 
\prod_{1\leq j<k\leq n} \frac{(g+\hat{\rho}_j+\hat{\rho}_k)_{\lambda_j+\lambda_k} (g+\hat{\rho}_j-\hat{\rho}_k)_{\lambda_j-\lambda_k}}{(\hat{\rho}_j+\hat{\rho}_k)_{\lambda_j+\lambda_k}(\hat{\rho}_j-\hat{\rho}_k)_{\lambda_j-\lambda_k}}  ,
\end{align*}
\begin{align*}
&V^{J,n}_{J_+,J_-} (\lambda ;g,g_0,g_1):=\prod_{j\in J} \frac{  (\hat{g}_0+\epsilon_j(\hat{\rho}_j+\lambda_j), \hat{g}_1+\epsilon_j(\hat{\rho}_j+\lambda_j))}{2\epsilon_j(\hat{\rho}_j+\lambda_j)(1+2\epsilon_j(\hat{\rho}_j+\lambda_j))} \\
&\times \prod_{\substack{j,j^\prime\in J\\ j<j^\prime}}
\frac{(g+\epsilon_j(\hat{\rho}_j+\lambda_j)+\epsilon_{j^\prime}(\hat{\rho}_{j^\prime
}+\lambda_{j^\prime}))(1+g+\epsilon_j(\hat{\rho}_j+\lambda_j)+\epsilon_{j^\prime}(\hat{\rho}_{j^\prime
}+\lambda_{j^\prime}))}{(\epsilon_j(\hat{\rho}_j+\lambda_j)+\epsilon_{j^\prime}(\hat{\rho}_{j^\prime
}+\lambda_{j^\prime}))(1+\epsilon_j(\hat{\rho}_j+\lambda_j)+\epsilon_{j^\prime}(\hat{\rho}_{j^\prime
}+\lambda_{j^\prime}))} \\
&\times \prod_{j\in J,k\in J^c} \frac{(g+\epsilon_j(\hat{\rho}_j+\lambda_j)+\hat{\rho}_k +\lambda_k)(g+\epsilon_j(\hat{\rho}_j+\lambda_j)-\hat{\rho}_k -\lambda_k)}{(\epsilon_j(\hat{\rho}_j+\lambda_j)+\hat{\rho}_k +\lambda_k)(\epsilon_j(\hat{\rho}_j+\lambda_j)-\hat{\rho}_k -\lambda_k)}
\end{align*}
with $\epsilon_j\equiv \epsilon_j(J_+,J_-)$,
\begin{align*}
 &U^{J,n}_{K,p} (\lambda ;g,g_0,g_1):=  (-1)^p \sum_{\substack{I_+,I_-\subset K\\ I_+\cap I_-=\emptyset \\  |I_+|+|I_-|= p }} 
  \Biggl( \prod_{j\in I}  \frac{   (\hat{g}_0+\epsilon_j(\hat{\rho}_j+\lambda_j),\hat{g}_1+\epsilon_j(\hat{\rho}_j+\lambda_j))}{2\epsilon_j(\hat{\rho}_j+\lambda_j)(1+2\epsilon_j(\hat{\rho}_j+\lambda_j))} \\
&\times \prod_{\substack{j,j^\prime\in I\\ j<j^\prime}}
\frac{(g+\epsilon_j(\hat{\rho}_j+\lambda_j)+\epsilon_{j^\prime}(\hat{\rho}_{j^\prime
}+\lambda_{j^\prime}))(1-g+\epsilon_j(\hat{\rho}_j+\lambda_j)+\epsilon_{j^\prime}(\hat{\rho}_{j^\prime
}+\lambda_{j^\prime}))}{(\epsilon_j(\hat{\rho}_j+\lambda_j)+\epsilon_{j^\prime}(\hat{\rho}_{j^\prime
}+\lambda_{j^\prime}))(1+\epsilon_j(\hat{\rho}_j+\lambda_j)+\epsilon_{j^\prime}(\hat{\rho}_{j^\prime
}+\lambda_{j^\prime}))}  \\
&\times \prod_{j\in I,k\in K\setminus I} \frac{(g+\epsilon_j(\hat{\rho}_j+\lambda_j)+\hat{\rho}_k +\lambda_k)(g+\epsilon_j(\hat{\rho}_j+\lambda_j)-\hat{\rho}_k -\lambda_k)}{(\epsilon_j(\hat{\rho}_j+\lambda_j)+\hat{\rho}_k +\lambda_k)(\epsilon_j(\hat{\rho}_j+\lambda_j)-\hat{\rho}_k -\lambda_k)}\Biggr)
\end{align*}
with $\epsilon_j\equiv \epsilon_j(I_+,I_-)$, $I=I_+\cup I_-$ and $p=0,\ldots ,|K|$.

\subsection{Branching polynomials}
Upon picking parameters of the form
\begin{subequations}
\begin{equation}\label{AW->Ja}
 t=q^{g}, \qquad  t_l=(-1)^l q^{h_l}\quad
(l=0,1,2,3)
\end{equation}
such that $h_0+h_2=g_0$ and $h_1+h_3=g_1$,
the Macdonald-Koornwinder polynomials degenerate into the symmetric Jacobi polynomials for $q\to 1$ \cite[\S 11]{mac:orthogonal}, \cite[Sec. 4.1]{die:commuting}, \cite[Sec. 4.3]{die:properties}:
\begin{equation}\label{AW->Jb}
P_\lambda^{J}(x_1,\ldots ,x_n; g,g_0,g_1)=\lim_{q\to 1} 
P_\lambda (x_1,\ldots , x_n; q,t,t_l) .
\end{equation}
\end{subequations}

The corresponding branching polynomials $P_{\lambda/\mu}(x; q,t,t_l)$ 
\eqref{br-pol}, \eqref{br-coef} are seen to
degenerate in this limit into the following branching polynomials for the symmetric Jacobi polynomials:
\begin{subequations}
\begin{equation}\label{br-pol-J}
\boxed{P^{J}_{\lambda/\mu}(x; g,g_0,g_1) =
\sum_{0\leq k\leq d} B_{\lambda/\mu}^{J,k} (g,g_0, g_1)
\sin^{2k}\left(\frac{x}{2}\right) }
\end{equation}
with
\begin{equation}\label{br-coef-J}
\boxed{B_{\lambda/\mu}^{J,k} (g,g_0,g_1)=
4^m (-1)^{|\lambda|-|\mu|}
C^{J,(n+1)^m-\lambda^\prime,m}_{n^m-\mu^\prime,m-k}({\textstyle \frac{1}{g},\frac{g_0}{g}, \frac{g_1}{g}}) }
 \end{equation}
 \end{subequations}
for $k=0,\ldots ,d=d(\lambda,\mu)$. Here one uses that
$$
\lim_{q\to 1} (-4)^{-k} \langle x, t_0\rangle_{q,k}=   \sin^{2k}\left(\frac{x}{2}\right) $$
and that
$$
\lim_{q\to 1}  4^{-r} C^{\mu,n}_{\lambda,r}(q,t,t_l)
=C^{J,\mu,n}_{\lambda,r}(g,g_0,g_1) .
$$

\section{Laguerre level}\label{L:sec}

\subsection{Symmetric Laguerre polynomials \cite{mac:hypergeometric,las:laguerre,bak-for:calogero-sutherland,die:confluent,xu:orthogonal,ari-dav-ola:differential,dum-ede-shu:MOPS,hal-lan:unified}}

The symmetric Laguerre polynomials are even polynomials
\begin{equation}
P_\lambda(x_1,\ldots ,x_n)=P_\lambda^{L}(x_1,\ldots ,x_n; g,h,\omega )\qquad ( \lambda\in \Lambda_n)
\end{equation}
of the form in Eqs.
\eqref{pa}, \eqref{pb}, with
$$M_\lambda(x_1,\ldots ,x_n)= m_\lambda (x_1^2,\ldots ,x_n^2) $$
and 
\begin{align}
\Delta(x_1,\ldots ,x_n) &=\Delta^{L}(x_1,\ldots ,x_n; g,h,\omega ) \\
&:=\prod_{1\leq j\leq n} e^{-\omega  x_j^2} \left| x_j
\right|^{2h-1}
\prod_{1\leq j<k\leq n} \left|  x_j^2-x_k^2
\right|^{2g}  \nonumber
\end{align}
supported on $\mathcal{D}=\mathbb{R}^n$. Here it is
assumed that the parameters $g, h$ and the scaling parameter $\omega$ are all positive.

\subsection{Pieri coefficients \cite{die:confluent}}

The Pieri coefficients $C^{\mu, n}_{\lambda , r}$ \eqref{pieri} for the symmetric Laguerre polynomials
 associated with the multiplication of $P^{L}_\lambda(x_1,\ldots ,x_n;g,h,\omega )$ by
\begin{equation}
E_r(x_1,\ldots ,x_n)=E_r^{L}(x_1,\ldots ,x_n)
:=\boldsymbol{e}_r(x_1^2 ,\ldots , x_n^2) 
\end{equation}
are given by 
\begin{equation}
C^{\mu, n}_{\lambda , r}=C^{L,\mu,n}_{\lambda , r}(g,h,\omega ):=  (-\omega )^{-r}
\frac{p^{L}_\lambda(g,h,\omega )}{p^{L}_\mu (g,h,\omega )} V^{L,n}_{J_+,J_-} (\lambda ;g,h) U^{L,n}_{J^c,r-|J|} (\lambda ; g,h) ,
\end{equation}
where
\begin{equation*}
p^{L}_\lambda(g,h,\omega ):= (-\omega)^{-|\lambda |}
 \prod_{1\leq j\leq n}  ((n-j)g + h)_{\lambda_j}
\prod_{1\leq j<k\leq n}  \frac{(1+(k-j)g)_{\lambda_j-\lambda_k}}{((k-j)g)_{\lambda_j-\lambda_k}} ,
\end{equation*}
\begin{align*}
&V^{L,n}_{J_+,J_-} (\lambda ;g,h):=  \prod_{j\in J_+}  ((n-j)g+h +\lambda_j)\prod_{j\in J_-}((n-j)g+\lambda_j)\\
\times 
&\prod_{\substack{j\in J_+\\ j'\in J_-}} 
\biggl(1+  \frac{g}{(j'-j)g +\lambda_j-\lambda_{j'}} \biggr) \biggl( 1+  \frac{g}{(j'-j)g +\lambda_j-\lambda_{j'}+1}\
\biggr )\\
\times 
&\prod_{\substack{j\in J_+\\ k\not\in J_+\cup J_-}} \biggl(1+  \frac{g}{(k-j)g +\lambda_j-\lambda_k})\biggr)
\prod_{\substack{j\in J_-\\k\not\in J_+\cup J_-}}  \biggl(1-  \frac{g}{(k-j)g +\lambda_j-\lambda_k} 
\biggr)
 ,
\end{align*}
and
\begin{align*}
 &U^{L,n}_{K,p} (\lambda ;g,h):=  (-1)^p \sum_{\substack{I_+, I_-\subset K\\ I_+\cap I_- = \emptyset \\  |I_+|+|I_-|=p}}
  \prod_{j\in I_+}  ((n-j)g+h +\lambda_j)\prod_{j\in I_-}((n-j)g+\lambda_j)\\
\times 
&\prod_{j\in I_+,j^\prime \in I_-} 
\biggl(1+  \frac{g}{(j'-j)g +\lambda_j-\lambda_{j'}} \biggr) \biggl( 1-  \frac{g}{(j'-j)g +\lambda_j-\lambda_{j'}+1}\
\biggr )\\
\times 
&\prod_{\substack{j\in I_+\\ k\in K\setminus (I_+\cup I_-)}} \biggl(1+  \frac{g}{(k-j)g+\lambda_j-\lambda_k})\biggr) 
\prod_{\substack{j\in I_-\\ k\in K\setminus (I_+\cup I_-)}} \biggl(1-  \frac{g}{(k-j)g +\lambda_j-\lambda_k} 
\biggr)
\end{align*}
for $p=0,\ldots ,|K|$.

\subsection{Branching polynomials}

Upon picking parameters such that
\begin{subequations}
\begin{equation}\label{W->La}
 g_0+g_1 =h\quad\text{and}\quad g_{l+2}=\frac{1}{ \omega_l\beta^2} \quad (l=0,1),
 \end{equation}
 with $ \omega_0+\omega_1=\omega$,
the symmetric Wilson polynomials degenerate into the symmetric Laguerre polynomials via the following limit
\cite[Sec. 4.1]{die:confluent}:
\begin{equation}\label{W->Lb}
P_\lambda^{L}(x_1,\ldots ,x_n; g,h,\omega )=\lim_{\beta \to 0}  \beta^{2|\lambda |}
P_\lambda^W ({\textstyle \frac{x_1}{\beta},\ldots , \frac{x_n}{\beta}; g,g_l}) .
\end{equation}
\end{subequations}

The corresponding branching polynomials $\beta^{2(|\lambda |-|\mu|)} P^W_{\lambda/\mu}(\frac{x}{\beta}; g,g_l )$ 
\eqref{br-pol-W}, \eqref{br-coef-W} are seen to
converge correspondingly to the following branching polynomials for the symmetric Laguerre polynomials:
\begin{subequations}
\begin{equation}\label{br-pol-L}
\boxed{
P^{L}_{\lambda/\mu}(x; g,h,\omega ) =
\sum_{0\leq k\leq d} B_{\lambda/\mu}^{L,k} (g,h,\omega )
x^{2k}}
\end{equation}
with
\begin{equation}\label{br-coef-L}
\boxed{B_{\lambda/\mu}^{L,k} (g,h,\omega )=
(-1)^{k+|\lambda|-|\mu|}  C^{L,(n+1)^m-\lambda^\prime,m}_{n^m-\mu^\prime,m-k}({\textstyle \frac{1}{g},\frac{h}{g},\frac{\omega}{g} }) }
 \end{equation}
 \end{subequations}
for $k=0,\ldots ,d=d(\lambda,\mu)$. Here one uses that
$$
\lim_{\beta\to 0} \beta^{2k}(g_0+ix\beta^{-1})_k(g_0-ix\beta^{-1})_k=  x^{2k}$$
and that
$$
\lim_{\beta \to 0} (-1)^r  \beta^{2(|\lambda|-|\mu|+r)}C^{W,\mu,n}_{\lambda,r}(g,g_l)
=C^{L,\mu,n}_{\lambda,r}(g,h,\omega ) .
$$

\section{Hermite level}\label{H:sec}

\subsection{Symmetric Hermite polynomials \cite{mac:hypergeometric,las:hermite,bak-for:calogero-sutherland,die:confluent,xu:orthogonal,dum-ede-shu:MOPS,hal-lan:unified}}

The symmetric Hermite polynomials
\begin{equation}
P_\lambda(x_1,\ldots ,x_n)=P_\lambda^{H}(x_1,\ldots ,x_n; g,\omega )\qquad ( \lambda\in \Lambda_n)
\end{equation}
are of the form in Eqs.
\eqref{pa}, \eqref{pb}, with
$$M_\lambda(x_1,\ldots ,x_n)= m_\lambda (x_1,\ldots ,x_n) $$
and 
\begin{align}
\Delta(x_1,\ldots ,x_n) &=\Delta^{H}(x_1,\ldots ,x_n; g,\omega )  \\
&:=\prod_{1\leq j\leq n} e^{-\omega x_j^2}
\prod_{1\leq j<k\leq n} \left|  x_j-x_k
\right|^{2g}  \nonumber
\end{align}
supported on $\mathcal{D}=\mathbb{R}^n$. Here it is
assumed that the parameter $g$ and the scaling parameter $\omega$ are both positive.

\subsection{Pieri coefficients \cite{die:confluent}}
Let us denote by $C^{H,\mu,n}_{\lambda , r}(g,\omega )$
the Pieri coefficients $C^{\mu, n}_{\lambda , r}$ \eqref{pieri}  associated with the multiplication of
the symmetric Hermite polynomial
$P^{H}_\lambda(x_1,\ldots ,x_n;g,\omega )$ by the elementary symmetric polynomial
\begin{equation}
E_r(x_1,\ldots ,x_n)=E_r^{H}(x_1,\ldots ,x_n)
:= \boldsymbol{e}_r(x_1 ,\ldots , x_n) .
\end{equation}
It is known that $C^{H,\mu,n}_{\lambda , r}(g,\omega )=0$ unless $\mu\sim_r\lambda$, but
unlike in the preceding cases above, a general explicit expression for this Pieri coefficient is not available
except when the cardinality of  $J=J(\lambda ,\mu)$ is equal to $r$:
\begin{equation}
C^{H,\mu,n}_{\lambda , r}(g,\omega )= 
\frac{p^{H}_\lambda(g)}{p^{H}_\mu (g)} V^{H,n}_{J_+,J_-} (\lambda ;g,\omega )  \qquad (\text{when}\  |J_+|+|J_-|=r),
\end{equation}
where
\begin{equation*}
p^{H}_\lambda(g):=  \prod_{1\leq j  < k\leq n} \frac{ ( 1+(k-j)g)_{\lambda_j - \lambda_k} }{ ((k-j)g)_{\lambda_j - \lambda_k}   }
\end{equation*}
and
\begin{align*}
&V^{H,n}_{J_+,J_-} (\lambda ;g,\omega ):=  
\prod_{j\in J_-}  \frac{(n-j)g+\lambda_j}{2\omega}\\
\times 
&\prod_{j\in J_+,j'\in J_-} 
\biggl(1+  \frac{g}{(j'-j)g +\lambda_j-\lambda_{j'}} \biggr) \biggl( 1+  \frac{g}{(j'-j)g+\lambda_j-\lambda_{j'}+1}\
\biggr )\\
\times 
&\prod_{\substack{j\in J_+\\ k\not\in J_+\cup J_-}} \biggl(1+  \frac{g}{(k-j)g+\lambda_j-\lambda_k}\biggr)
\prod_{\substack{j\in J_-\\ k\not\in J_+\cup J_-}}  \biggl(1-  \frac{g}{(k-j)g+\lambda_j-\lambda_k} 
\biggr) .
\end{align*}

\subsection{Branching polynomials}
Upon picking parameters such that
\begin{subequations}
\begin{equation}\label{cH->Ha}
 g_l= \frac{1}{\omega_l \beta^2} \quad (l=0,1)
\end{equation}
with $\omega_0+\omega_1=\omega$, the symmetric continuous Hahn polynomials degenerate into the symmetric Hermite polynomials via the following limit
\cite[Sec. 4.1]{die:confluent}:
\begin{equation}\label{cH->Hb}
P_\lambda^{H}(x_1,\ldots ,x_n; g,\omega )=\lim_{\beta \to 0}  \beta^{|\lambda |}
P_\lambda^{cH} ({\textstyle \frac{x_1}{\beta},\ldots , \frac{x_n}{\beta}; g,g_0,g_1}) .
\end{equation}
\end{subequations}

Since it turns out to be very cumbersome to perform the limit \eqref{cH->Ha}, \eqref{cH->Hb} at the level of the Pieri coefficients, we have done so instead at the level of the Cauchy identity (see Appendix \ref{cauchy:app}).  
This allows to deduce the branching formula for the symmetric Hermite polynomials
directly following the approach in \cite[Sec. 4]{die-ems:branching} for the Macdonald-Koornwinder polynomials.

Specifically, by expanding the first two factors of the trivial identity
$$
 \prod_{\substack{1\leq j\leq m\\1\leq k\leq n+1}}  (x_j-z_k) =
  \prod_{\substack{1\leq j\leq m\\1\leq k\leq n}}  (x_j-z_k)  \prod_{{1\leq j\leq m}}  (x_j-z_{n+1}) 
$$
with the aid of the Cauchy identity for the symmetric Hermite polynomials in Eq. \eqref{H:cauchy}, and employing the elementary expansion
\begin{equation}\label{r-c:cauchy}
\prod_{{1\leq j\leq m}}  (x_j-z) =\sum_{0\leq r\leq m} (-1)^{m-r}\boldsymbol{e}_r(x_1,\ldots ,x_m) z^{m-r},\qquad z=z_{n+1},
\end{equation}
for the last factor, one arrives at the equality
\begin{equation*}
\sum_{\lambda\subset (n+1)^m} (-1)^{m(n+1)-|\lambda |}
P_\lambda^H (x_1,\ldots ,x_m; g,\omega ) P_{m^{n+1}-\lambda^\prime}^H (z_1,\ldots ,z_{n+1};{\textstyle \frac{1}{g},\frac{\omega}{g}}) 
\end{equation*}
\begin{align*}
 = \sum_{\substack{\mu\subset n^m\\ 0\leq r\leq m}} (-1)^{m(n+1)-|\mu |-r}& \Bigl( 
P_{m^n-\mu^\prime}^H (z_1,\ldots ,z_n; {\textstyle \frac{1}{g},\frac{\omega}{g}}) z_{n+1}^{m-r}  \\
& \times  \boldsymbol{e}_{r}(x_1,\ldots,x_m) P_{\mu}^H (x_1,\ldots ,x_m;g,\omega )
  \Bigr) .
\end{align*}
After invoking the Pieri formula for the symmetric Hermite polynomials and reordering of the summations, the RHS is rewritten as
\begin{align*}
 =  & \sum_{\lambda\subset (n+1)^m }  (-1)^{m(n+1)-|\lambda |}  \Bigl( P_{\lambda}^H (x_1,\ldots ,x_m;g,\omega )  \\
& \times \sum_{\substack{\mu\subset n^m,  \, 0\leq r\leq m\\ \mu\sim_{r} \lambda}} 
 (-1)^{r+|\lambda|-|\mu|} C^{H,\lambda,m}_{\mu , r}(g,\omega )   
P_{m^n-\mu^\prime}^H (z_1,\ldots ,z_n;    {\textstyle \frac{1}{g},\frac{\omega}{g}})  z_{n+1}^{m-r}    \Bigr)  .
\end{align*}
Hence, it is seen by comparing with the LHS that for any $\lambda\subset (n+1)^m$:
\begin{align*}
&P_{m^{n+1}-\lambda^\prime}^H (z_1,\ldots ,z_{n+1};{\textstyle \frac{1}{g},\frac{\omega}{g}} ) 
 =    \\
&  \sum_{\substack{\mu\subset n^m, \, 0\leq r\leq m\\ \mu\sim_{r} \lambda }} 
 (-1)^{r+|\lambda|-|\mu|} C^{H,\lambda,m}_{\mu , r}(g,\omega )   
P_{m^n-\mu^\prime}^H (z_1,\ldots ,z_n;{\textstyle \frac{1}{g},\frac{\omega}{g}}) z_{n+1}^{m-r}      ,
\end{align*}
i.e. for any $\lambda\in m^{n+1}$ (cf. \cite[Lem. 5]{die-ems:branching}):
\begin{align*}
&P^H_{\lambda} (z_1,\ldots ,z_{n+1};g,\omega) 
 =    \\
&  \sum_{\substack{\mu\subset m^n,\, \mu\preceq\lambda \\ m-d(\lambda ,\mu) \leq r\leq m }} 
 (-1)^{m-r+|\lambda|-|\mu|} 
  C^{H, (n+1)^m-\lambda^\prime,m}_{n^m-\mu^\prime,r}({\textstyle \frac{1}{g},\frac{\omega}{g}}) 
P^H_{\mu} (z_1,\ldots ,z_n;g,\omega)  z_{n+1}^{m-r}    .  \nonumber
\end{align*}

The upshot is  that the branching rule for the symmetric Hermite polynomials
is of the form in Eq. \eqref{bf} with the
one-variable branching polynomial given by
\begin{subequations}
\begin{equation}
\boxed{
P^H_{\lambda/\mu}(x; g,\omega )=
\sum_{0\leq k\leq d} B_{\lambda /\mu}^{H,k} (g,\omega )
x^k}
\end{equation}
where
\begin{equation}
\boxed{
B_{\lambda/\mu}^{H, k} (g,\omega ) 
  =  (-1)^{ k + |\mu| - |\lambda| }  C^{H, (n+1)^m-\lambda^\prime,m}_{n^m-\mu^\prime,m-k}({\textstyle \frac{1}{g},\frac{\omega}{g}}) }
  \end{equation}
  \end{subequations}
for $k=0,\ldots ,d=d(\lambda,\mu)$.

\begin{remark}\label{br-proof:rem}
Given the pertinent Pieri coefficients, one may alternatively derive
the branching rules in this paper directly from the corresponding Cauchy identities by adapting the above proof for the Hermite case.
At the Askey-Wilson level, this boils down to replacing in the proof at issue:
(i) the degenerate Cauchy identity \eqref{H:cauchy} by Mimachi's Cauchy identity
\eqref{AW:cauchy} and (ii) the elementary 
expansion \eqref{r-c:cauchy} by
the special `column-row'  case of Okounkov's Cauchy identity in
 \cite[Lem. 5.1]{kom-nou-shi:kernel}:
\begin{align}\label{rc-cauchy:AW}
\prod_{1\leq j\leq m}&
( e^{ix_j}+e^{-ix_j}-e^{iz}-e^{-iz} )= \\
&\sum_{0\leq r\leq m} (-1)^{m-r} E_{r}(x_1,\ldots,x_m;t,t_0) \,  \langle z;t_0\rangle_{t,m-r} . \nonumber
\end{align}
Indeed, this way one precisely reproduces the proof of the branching rule for the Macdonald-Koornwinder polynomials in \cite[Sec. 4]{die-ems:branching}. The branching rules for the remaining hypergeometric families follow in turn with the aid of the degenerate Cauchy identities in Appendix \ref{cauchy:app} and the corresponding degenerations of the `column-row' Cauchy identity \eqref{rc-cauchy:AW}. 
At the Wilson level  and the continuous Hahn level the degenerations of the latter `column-row'
Cauchy identity become of the form
\begin{align}\label{rc-cauchy:W}
\prod_{1\leq j\leq m} & (x_j^2- z^2)=\\
&(-1)^{m} \sum_{0\leq r\leq m} E_{r}^W(x_1,\ldots,x_m;g,g_0) \, g^{2(m-r)}\,
 ({\textstyle \frac{g_0+iz}{g},\frac{g_0-iz}{g}})_{m-r} \nonumber
\end{align}
and
\begin{equation}\label{rc-cauchy:cH}
\prod_{1\leq j\leq m} (x_j - z)
=
i^m \sum_{0\leq r\leq m}  E^{cH}_r(x_1,\dots, x_m;g,g_0) \, g^{m-r} \, ({\textstyle \frac{g_0+iz}{g}})_{m-r},
\end{equation}
whereas at the Jacobi level and the Laguerre level the degenerate Cauchy identity is of the elementary form in \eqref{r-c:cauchy},  up to a trigonometric change of variables $x_j\to \sin^2 \left(\frac{x_j}{2}\right)$, $z\to\sin^2\left(\frac{z}{2}\right)$
and a quadratic change of variables $x_j\to x_j^2$, $z\to z^2$, respectively.
\end{remark}

\begin{remark}
It is known that the highest-degree leading terms of $P_\lambda^H (x_1,\ldots ,x_n;g,\omega)$ consist of the (monic) Jack polynomial
$P_\lambda (x_1,\ldots ,x_n;g)$ \cite{mac:hypergeometric,las:hermite,bak-for:calogero-sutherland,die:confluent}. By filtering the highest-degree terms on both sides of the branching formula \eqref{bf} for the symmetric Hermite polynomial $P_\lambda^H (x_1,\ldots ,x_n;g,\omega)$, one recovers in turn a celebrated branching rule for the Jack polynomials \cite{sta:some,mac:symmetric,oko-ols:shifted}:
$$P_\lambda (x_1,\ldots ,x_n,x;g)=
\sum_{\substack{\mu\in\Lambda_n,\,\mu\subset\lambda \\ \lambda/\mu\ \text{horizontal\ strip}}}  P_\mu (x_1,\ldots ,x_n;g)
P_{\lambda/\mu}(x;g)
$$
$(\lambda\in\Lambda_{n+1})$, where
$$
P_{\lambda/\mu}(x;g )= x^{|\lambda |-|\mu |} B_{\lambda/\mu} (g)
$$
and
\begin{align*}
B_{\lambda/\mu} (g)=&B^{H,d}_{\lambda/\mu}(g,\omega)=C^{H,(n+1)^m-\lambda^\prime,m}_{n^m-\mu^\prime ,m-d}({\textstyle \frac{1}{g},\frac{\omega}{g}}) 
\quad \text{with}\ d=|\lambda|-|\mu|
\\
= & \frac{p^H_{n^m-\mu^\prime}(\frac{1}{g})}{p^H_{(n+1)^m-\lambda^\prime}(\frac{1}{g})}
V^{H,m}_{J(n^m-\mu^\prime ,(n+1)^m-\lambda^\prime),\emptyset } (n^m-\mu^\prime ;{\textstyle \frac{1}{g},\frac{\omega}{g}})  \\
=&
\prod_{1\leq j<k\leq m} \frac{(1+(k-j)g^{-1})_{\mu_j^\prime-\mu_k^\prime}}{((k-j)g^{-1})_{\mu_j^\prime-\mu_k^\prime}} 
 \frac{((k-j)g^{-1})_{\lambda_j^\prime-\lambda_k^\prime}}{(1+(k-j)g^{-1})_{\lambda_j^\prime-\lambda_k^\prime}} \\
 &\times \prod_{\substack{ 1\leq j,k\leq m\\ \mu_j^\prime\neq \lambda_j^\prime \\ \mu_k^\prime = \lambda_k^\prime}}
 \biggl(1+  \frac{g^{-1}}{(k-j)g^{-1}+\mu_j^\prime-\mu_k^\prime }\biggr) \\
 =&
  \prod_{\substack{ 1\leq j<k\leq m\\ \mu_j^\prime =\lambda_j^\prime \\ \mu_k^\prime \neq \lambda_k^\prime}}
 \biggl(1+\frac{1}{k-j+g(\lambda_j^\prime-\lambda_k^\prime  )}\biggr) 
 \biggl(1-\frac{1}{k-j+g(\mu_j^\prime-\mu_k^\prime  )}\biggr)  \\
 =&
 \prod_{1\leq j\leq k\leq \ell (\mu)} \frac{(\mu_j-\mu_k+g(1+k-j),1+\mu_j-\lambda_{k+1}+g(k-j))_{\lambda_j-\mu_j}}{(1+\mu_j-\mu_k+g(k-j),\mu_j-\lambda_{k+1}+g(1+k-j))_{\lambda_j-\mu_j}} .
 \end{align*}
 This branching formula for the Jack polynomials amounts to the $t=q^g$, $q\to 1$ degeneration of the branching rule for the Macdonald polynomials, cf. Remark \ref{mac-br:rem} of Section \ref{AW:sec} above.
\end{remark}

\appendix

\section{Hypergeometric Cauchy identities}\label{cauchy:app}
This appendix collects the degenerations of Mimachi's (dual) Cauchy identity for all families of symmetric hypergeometric polynomials
considered above. At the lowest level of the symmetric Hermite polynomials, we relied on the pertinent Cauchy identity  for a direct verifcation of our branching rule in the absence of explicit limiting expressions for the corresponding Pieri coefficients.

\subsection{Askey-Wilson level \cite[Thm. 2.1]{mim:duality}}
\begin{align}\label{AW:cauchy}
& \prod_{\substack{1\leq j\leq m\\1\leq k\leq n}}  (e^{ix_j}+e^{-ix_j}-e^{iz_k}-e^{-iz_k}) = \\
& \sum_{\lambda\subset n^m} (-1)^{mn-|\lambda |}
P_\lambda (x_1,\ldots ,x_m;q,t,t_l) P_{m^n-\lambda^\prime}  (z_1,\ldots ,z_n;t,q,t_l) \nonumber
\end{align}
When $t=q$, $t_0=-t_1=q^{1/2}$ and $t_2=-t_3=q$, Mimachi's Cauchy identity \eqref{AW:cauchy} recovers a well-known Cauchy identity for the
symplectic Schur functions \cite{mor:spin,kin:branching,ter:robinson-schensted,ham-kin:bijective}.

\subsection{Wilson level}
\begin{align}\label{W:cauchy}
& \prod_{\substack{1\leq j\leq m\\1\leq k\leq n}}  (x_j^2-z_k^2) = \\
& \sum_{\lambda\subset n^m} (-g^2)^{mn-|\lambda |}
P_\lambda^W (x_1,\ldots ,x_m;g,g_l) P_{m^n-\lambda^\prime}^W  ({\textstyle \frac{z_1}{g},\ldots ,\frac{z_n}{g}; \frac{1}{g},\frac{g_l}{g}}) \nonumber
\end{align}
Eq. \eqref{W:cauchy} is obtained from Eq. \eqref{AW:cauchy} via the limit transition \eqref{AW->Wa}, \eqref{AW->Wb}.

\subsection{Continuous Hahn level}
\begin{align}\label{cH:cauchy}
& \prod_{\substack{1\leq j\leq m\\1\leq k\leq n}}  (x_j-z_k) = \\
& \sum_{\lambda\subset n^m} (-g)^{mn-|\lambda |}
P_\lambda^{cH} (x_1,\ldots ,x_m;g,g_0,g_1) P_{m^n-\lambda^\prime}^{cH}  ({\textstyle \frac{z_1}{g},\ldots ,\frac{z_n}{g};\frac{1}{g},\frac{g_0}{g},\frac{g_1}{g}}) \nonumber
\end{align}
Eq. \eqref{cH:cauchy} is obtained from Eq. \eqref{AW:cauchy} via the limit transition \eqref{AW->cHa}, \eqref{AW->cHb}.

\subsection{Jacobi level \cite[Sec. 6]{ser:some}, \cite[Thm. 4.1]{mim:duality}}
\begin{align}\label{J:cauchy}
&\prod_{\substack{1\leq j\leq m\\1\leq k\leq n}}  (e^{ix_j}+e^{-ix_j}-e^{iz_k}-e^{-iz_k}) = \\
& \sum_{\lambda\subset n^m} (-1)^{mn-|\lambda |}
P_\lambda^J (x_1,\ldots ,x_m;g,g_0,g_1) P_{m^n-\lambda^\prime}^J  (z_1,\ldots ,z_n;{\textstyle\frac{1}{g},\frac{g_0}{g},\frac{g_1}{g}}) \nonumber
\end{align}
Eq. \eqref{J:cauchy}  is obtained from Eq. \eqref{AW:cauchy} via the limit transition \eqref{AW->Ja}, \eqref{AW->Jb} \cite{mim:duality}.

\subsection{Laguerre level}
\begin{align}\label{L:cauchy}
& \prod_{\substack{1\leq j\leq m\\1\leq k\leq n}}  (x_j^2-z_k^2) = \\
& \sum_{\lambda\subset n^m} (-1)^{mn-|\lambda |}
P_\lambda^L (x_1,\ldots ,x_m;g,h,\omega ) P_{m^n-\lambda^\prime}^L (z_1,\ldots ,z_n;{\textstyle \frac{1}{g},\frac{h}{g},\frac{\omega}{g}}) \nonumber
\end{align}
Eq. \eqref{L:cauchy} is obtained from Eq. \eqref{W:cauchy} via the limit transition \eqref{W->La}, \eqref{W->Lb}.

\subsection{Hermite level}
\begin{align}\label{H:cauchy}
& \prod_{\substack{1\leq j\leq m\\1\leq k\leq n}}  (x_j-z_k) = \\
& \sum_{\lambda\subset n^m} (-1)^{mn-|\lambda |}
P_\lambda^H (x_1,\ldots ,x_m;g,\omega ) P_{m^n-\lambda^\prime}^H  (z_1,\ldots, z_n; {\textstyle \frac{1}{g}, \frac{\omega}{g}}) \nonumber
\end{align}
Eq. \eqref{H:cauchy} is obtained from Eq. \eqref{cH:cauchy} via the limit transition \eqref{cH->Ha}, \eqref{cH->Hb}.

\bibliographystyle{amsplain}

\end{document}